\documentclass[12pt,oneside,reqno]{amsart}
\usepackage{amsmath,latexsym,amssymb,amsthm,amsfonts,stmaryrd,mathrsfs,enumerate,graphicx}
\usepackage[colorlinks=true,linkcolor=blue,citecolor=blue,urlcolor=blue]{hyperref}
\newcommand{\be}{\begin{eqnarray}}
\newcommand{\ee}{\end{eqnarray}}
\newcommand{\ce}{\begin{eqnarray*}}
\newcommand{\de}{\end{eqnarray*}}
\newtheorem{theorem}{Theorem}[section]
\newtheorem{lemma}[theorem]{Lemma}
\newtheorem{remark}[theorem]{Remark}
\newtheorem{definition}[theorem]{Definition}
\newtheorem{proposition}[theorem]{Proposition}
\newtheorem{Examples}[theorem]{Examples}
\newtheorem{corollary}[theorem]{Corollary}

\def\[{{\Big[}}
\def\]{{\Big]}}
\def\<{{\langle}}
\def\>{{\rangle}}
\def\({{\Big(}}
\def\){{\Big)}}

\def\bt{\begin{theorem}}
\def\et{\end{theorem}}
\def\bl{\begin{lemma}}
\def\el{\end{lemma}}
\def\br{\begin{remark}}
\def\er{\end{remark}}
\def\bx{\begin{Examples}}
\def\ex{\end{Examples}}
\def\bd{\begin{definition}}
\def\ed{\end{definition}}
\def\bp{\begin{proposition}}
\def\ep{\end{proposition}}
\def\bc{\begin{corollary}}
\def\ec{\end{corollary}}

\def\cE{{\mathcal E}}

\def\mD{{\mathbb D}}

\def\mL{{\mathbb L}}

\def\mN{{\mathbb N}}

\def\mR{{\mathbb R}}

\def\mZ{{\mathbb Z}}
\def\geq{\geqslant}
\def\leq{\leqslant}

\allowdisplaybreaks

\begin{document}
\title{Boundedness of the Density of Wiener Functionals with Low Regularity}
\author{Hua Zhang}
\date{}
\dedicatory{School of Statistics and Data Science, Jiangxi University of Finance and Economics,\\
	Nanchang, Jiangxi 330013, P.R.China\\
	Email: H. Zhang: zh860801@163.com}
\thanks{This work is supported by National Natural Science Foundation of China (Grant Nos. 12261038 and 12371152), and Natural Science Foundation of Jiangxi Province (Grant Nos. 20232BAB201004 and 20242BAB23003).}

\begin{abstract}
In this paper we prove that if a Wiener functional $F$ is $1+\delta$-differentiable in the Malliavin sense for some $\delta>0$ and satisfies the nondegeneracy condition, then the density of the law of $F$ exists, is bounded and is H\"{o}lder continuous of every exponent $\beta<\delta$. This gives a partial positive answer to a problem proposed by Watanabe \cite{W}; to date, this seems to be the best result available for the boundedness of the density of Wiener functionals. The approach is based on the fractional order Sobolev spaces associated with a probability measure and on the Riesz transform.
\end{abstract}

\maketitle

\section{Introduction}
In recent years, the analysis of the density of Wiener functionals has drawn increasing attention, and the Malliavin calculus has proved to be an effective tool for dealing with this problem. For the existence of the density, Bouleau and Hirsch \cite{BH} established the absolute continuity criterion, which asserts that if $F\in\mD_1^2(\mR^d)$ and $\det\Sigma_F>0$ a.e., then the law of $F$ is absolutely continuous with respect to the Lebesgue measure. Using the canonical integration by parts formulas, one also knows that if $F\in\mD_{d+2}^{\infty-}(\mR^d)$ satisfies the nondegeneracy condition, then the density exists and is bounded; for more details the reader is referred to \cite{N} or \cite{MT}. More recently, Bally and Caramellino \cite{BC1}, using the Riesz transform introduced by Malliavin and Thalmaier \cite{MT}, proved that if $F\in\mD_2^{\infty-}$ satisfies the nondegeneracy condition, then the density $p_F$ exists and is bounded.

Watanabe \cite{W} has introduced the fractional order Sobolev spaces on Wiener space and proved the following theorem (see \cite[Theorem 2.2]{W}):
\bt
Let $F: B\rightarrow\mR^d$ be a Wiener map such that
\begin{enumerate}
\item[(i)] $F\in\mD_{1+\delta}^{\infty-}(\mR^d)$ for $\delta>0$;
\item[(ii)] $(\det(\Sigma_F))^{-1}\in\mL^{\infty-}$;
\item[(iii)] the density $p_F$ of the law of $F$ is bounded.
\end{enumerate}
Then, for every $\alpha>0$, every $0\leq\beta<\delta\wedge\alpha$ and $1<p<\infty$ satisfying
\ce
d-p\{d-(\alpha\wedge\delta-\beta)\}>0,
\de
we have $(1-\Delta)^{\beta/2}\delta_x\circ F\in\mD_p^{-\alpha}$ boundedly and continuously in $x$.
\et
However, assumption (iii) in the above theorem is redundant in two important regimes: if $1+\delta\geq d+2$, it follows from assumptions (i) and (ii) by the classical integration by parts argument, and by the result of Bally and Caramellino \cite{BC1} it is in fact implied by (i) and (ii) whenever $1+\delta\geq2$. We are therefore led to ask whether, for every $\delta>0$, assumption (iii) follows from assumptions (i) and (ii).

In the present paper we give a partial positive answer to this question. The conclusion is striking: for a nondegenerate $F\in\mD_1^{\infty-}$, an arbitrarily small additional amount of regularity suffices to obtain boundedness of the density of $F$.

Our main result can now be stated as follows. We shall use the following notation. $Q_d$ denotes the Poisson kernel on $\mR^d$, that is, the fundamental solution of $\Delta Q_d=\delta_0$ (the explicit formula is recalled in Section 2), and $H_i(F,G)$, $i=1,\cdots,d$, denotes the weight appearing in the integration by parts formula
\ce
E[(\partial_if)(F)G]=E[f(F)H_i(F,G)],\qquad f\in C_c^1(\mR^d),
\de
which will be introduced in Subsection 3.2 below. For a nondegenerate Wiener functional $F$, one has $H_i(F,G)=\delta(\gamma_F^{ij}G\,DF_j)$, where $\gamma_F^{ij}$ are the entries of the inverse Malliavin covariance matrix $\Gamma_F=\Sigma_F^{-1}$ (see Section 2) and $\delta$ is the Skorohod divergence operator. Throughout the paper the Einstein summation convention is used for repeated indices, and any norm involving a free index is understood as the maximum (equivalently, the sum) over all its values.

\bt\label{thm:main}
\begin{enumerate}
  \item Let $0<\delta\leq 1$. Let $F:B\rightarrow\mR^d$ be a Wiener functional such that $F\in\mD_{1+\delta}^{\infty-}(\mR^d)$. Suppose further that $F$ is nondegenerate, i.e.,
      \ce
      (\det(\Sigma_F))^{-1}\in\mL^{\infty-}.
      \de
Then the law of $F$ is absolutely continuous with respect to the Lebesgue measure on $\mR^d$ and the density can be represented as
\ce
p_F(x)=E[\sum_{i=1}^d\partial_iQ_d(F-x)H_i(F,1)].
\de
Moreover, $p_F$ is H\"{o}lder continuous of every exponent $\beta<\delta$, and for every $p>1$ there exists a positive constant $C=C(\delta,p,d,\|\gamma_F^{ij}DF_j\|_{\delta,p})$ such that
\be\label{w1}
\|p_F\|_{\infty}\leq C(\delta,p,d,\|\gamma_F^{ij}DF_j\|_{\delta,p})\|\gamma_F^{ij}DF_j\|_{\delta,p}.
\ee
  \item Let $0<\delta\leq 1$, $m\geq2$. Let $F:B\rightarrow\mR^d$ be a Wiener functional such that $F\in\mD_{m+\delta}^{\infty-}(\mR^d)$. Suppose further that $F$ is nondegenerate.
Then $p_F\in C^{m-1}(\mR^d)$ and for every multiindex $\alpha$ with $|\alpha|\leq m$ one has
\ce
\partial_{\alpha}p_F(x)=(-1)^{|\alpha|}E[\sum_{i=1}^d\partial_iQ_d(F-x)H_{(\alpha,i)}(F,1)].
\de
Moreover, for $|\alpha|\leq m-1$, $\partial_{\alpha}p_F$ is H\"{o}lder continuous of any exponent $\beta<\delta$ and for every $p>1$ there exists a positive constant $C=C(\delta,p,d,\|\gamma_F^{ij}DF_j\|_{m-1+\delta,p})$ such that
\be\label{w2}
\|\partial_{\alpha}p_F\|_{\infty}\leq C(\delta,p,d,\|\gamma_F^{ij}DF_j\|_{m-1+\delta,p})\|\gamma_F^{ij}DF_j\|_{m-1+\delta,p}.
\ee
\end{enumerate}
\et

\begin{remark}\label{rem:delta-sharp}
	The condition $\delta>0$ in Theorem~\ref{thm:main} is essential. The assumption
	$F\in\mD_1^{\infty-}$ together with the nondegeneracy condition is not
	sufficient to guarantee that the density $p_F$ is bounded.
	
	Let $W_1$ be a standard Gaussian random variable on an abstract Wiener space.
	Define a continuous function $\psi:\mathbb{R}\to\mathbb{R}$ by
	\ce
	\psi(x)=
	\begin{cases}
		\displaystyle \frac{1}{\log(1/|x|)}, & 0<|x|<e^{-1},\\[6pt]
		1, & |x|\ge e^{-1},\\[4pt]
		0, & x=0,
	\end{cases}
	\de
	and set
	\ce
	\phi(x)=\int_0^x \psi(t)\,dt,\qquad F=\phi(W_1).
	\de
	Then $\phi\in C^1(\mathbb{R})$ and $\phi'(x)=\psi(x)$.
	
	We verify the following facts.
	
	\begin{enumerate}
		\item $F\in\mD_1^{\infty-}$. Since $0\leq\psi\leq1$ and $\phi(0)=0$, we have $|\phi(x)|\leq|x|$, hence $F\in\mL^p$ for every $p>1$. Moreover,
		\ce
		DF=\psi(W_1)DW_1,
		\de
		and since $\psi$ is bounded, $DF\in L^p(H)$ for every $p>1$.
		Hence $F\in\mathbb{D}_1^p$ for all $p>1$, i.e., $F\in\mD_1^{\infty-}$.
		
		\item $F$ is nondegenerate. In the one-dimensional case,
		\ce
		\Sigma_F=|\psi(W_1)|^2,\qquad
		\Gamma_F=\Sigma_F^{-1}=|\psi(W_1)|^{-2}.
		\de
		For $|x|<e^{-1}$,
		\ce
		|\psi(x)|^{-2}=\log^2(1/|x|),
		\de
		so that
		\ce
		E[\Gamma_F^p]
		\le C+C\int_0^{e^{-1}}\log^{2p}(1/x)\,dx<\infty
		\de
		for every $p>1$. Thus $\Gamma_F\in\bigcap_{p>1}L^p$, which means
		that $F$ is nondegenerate.
		
		\item However, the density $p_F$ is unbounded. Since $\phi$ is strictly
		increasing,
		\ce
		p_F(y)=\frac{p_{W_1}(\phi^{-1}(y))}{\phi'(\phi^{-1}(y))}
		=\frac{p_{W_1}(\phi^{-1}(y))}{\psi(\phi^{-1}(y))}.
		\de
		As $y\to0$, we have $\phi^{-1}(y)\to0$, hence
		$\psi(\phi^{-1}(y))\to0$, while $p_{W_1}(0)>0$. Therefore
		\ce
		p_F(y)\to+\infty\qquad (y\to0),
		\de
		and consequently $p_F\notin L^\infty(\mathbb{R})$.
	\end{enumerate}
	
	Thus, merely assuming $F\in\mD_1^{\infty-}$ and nondegeneracy does not
	ensure that the density is bounded. This shows that the requirement
	$F\in\mD_{1+\delta}^{\infty-}$ with $\delta>0$ in
	Theorem~\ref{thm:main} cannot be weakened to $F\in\mD_1^{\infty-}$.
	
	Finally, we note that if one drops the nondegeneracy assumption, then
	even the existence of a density may fail. For instance, let
	\ce
	F=W_1+|W_1|.
	\de
	Then $F\in\mD_1^{\infty-}$, but
	\ce
	P(F=0)=P(W_1<0)=\frac12,
	\de
	so the law of $F$ has an atom and no density exists. In this case
	\ce
	\Sigma_F=|1+\operatorname{sgn}(W_1)|^2
	\de
	vanishes on $\{W_1<0\}$, so the nondegeneracy condition fails.
\end{remark}

The paper is organized as follows. In Section 2 we recall some basic facts about the Malliavin calculus and the fractional order Sobolev spaces on Wiener space. Section 3 is devoted to the abstract functional analytic framework. In Subsection 3.1, after introducing the fractional order Sobolev spaces associated with a probability measure, we first state our main abstract results on the representation and regularity of densities (Theorems \ref{elementary} and \ref{regularity}, together with Corollary \ref{holderdensity}); we then establish a sequence of auxiliary lemmas based on the Riesz transform, including the key boundedness estimate (Lemma \ref{boundmain}), and prove the stated results. Subsection 3.2 introduces the analogous Sobolev spaces associated with a Wiener functional and states the corresponding density estimates (Theorems \ref{boundness21} and \ref{boundness22}). Section 4 is devoted to the proof of Theorem \ref{thm:main}.

\section{Preliminaries on the Malliavin calculus}
Let $(B,H,\mu)$ be an abstract Wiener space.
We refer to \cite{N} (see also \cite{HY} and \cite{PM}) for the background on the Malliavin calculus. We denote the
gradient operator (or Shigekawa's $H$-derivative) by $D$, its dual divergence operator (or the Skorohod operator) by $\delta$ and
the Ornstein-Uhlenbeck operator by $L:=-\delta D$. The Sobolev spaces ${\mathbb D}_{\alpha}^p$, $1<p<\infty$, $\alpha\in{\mathbb R}$,
of real Wiener and generalized Wiener functionals are defined by
\begin{eqnarray*}
{\mathbb D}_{\alpha}^{p}=(1-L)^{-\alpha/2}(\mathbb{L}^p)
\end{eqnarray*}
with the norm
\begin{eqnarray*}
\|F\|_{\alpha,p}=\|(1-L)^{\alpha/2}F\|_{\mathbb{L}^p},
\end{eqnarray*}
where $\mathbb{L}^p$ is the usual $\mathbb{L}^p$-space. In particular
\begin{eqnarray*}
{\mathbb D}_0^p=\mathbb{L}^p,\quad{\mathbb D}_{\alpha'}^{p'}\subseteq{\mathbb D}_{\alpha}^{p}\quad\text{if}\quad p<p'\quad\text{and}
\quad\alpha<\alpha',\quad({\mathbb D}_{\alpha}^p)'={\mathbb D}_{-\alpha}^{p/(p-1)}.
\end{eqnarray*}
As in the Schwartz theory of distributions, we define the space of test Wiener functionals by
\begin{eqnarray*}
{\mathbb D}_{\infty}^{\infty-}=\bigcap_{\alpha>0}\bigcap_{1<p<\infty}{\mathbb D}_{\alpha}^p.
\end{eqnarray*}
Its dual, the space of generalized Wiener functionals, is given by
\begin{eqnarray*}
{\mathbb D}_{-\infty}^{1+}=\bigcup_{\alpha>0}\bigcup_{1<p<\infty}{\mathbb D}_{-\alpha}^p.
\end{eqnarray*}
For $\alpha\in{\mathbb R}$ we also use the notation
\begin{eqnarray*}
{\mathbb D}_{\alpha}^{\infty-}=\bigcap_{1<p<\infty}{\mathbb D}_{\alpha}^p,\qquad
{\mathbb L}^{\infty-}=\bigcap_{1<p<\infty}{\mathbb L}^p.
\end{eqnarray*}
If we consider, more generally, $E$-valued functionals, $E$ being a real separable Hilbert space, the corresponding spaces are
denoted by $\mathbb{L}^p(E)$, ${\mathbb D}_s^p(E)$, etc. In the sequel we shall use the following convention:
\begin{eqnarray*}
\text{for}\quad F\in{\mathbb D}_{-\infty}^{1+},\quad\text{we set}\quad\|F\|_{\alpha,p}=\infty\quad\text{if}\quad F
\notin{\mathbb D}_{\alpha}^p.
\end{eqnarray*}

We will need the following result (see Kusuoka \cite{K}). For the spaces ${\mathbb D}_{\alpha}^p$,
multiplication is continuous: for every $\alpha\geqslant 0$ and $1<p,q,r<\infty$ such that $1/p+1/q=1/r$,
the map $(F,G)\in {\mathbb D}_{\alpha}^p\times{\mathbb D}_{\alpha}^q\rightarrow F\cdot G\in{\mathbb D}_{\alpha}^r$
is continuous, i.e., there exists a positive constant $C=C(p,q;\alpha)$ such that
\begin{eqnarray}
\|F\cdot G\|_{\alpha,r}\leqslant C\|F\|_{\alpha,p}\cdot\|G\|_{\alpha,q}\quad\text{for all}\quad F,G\in
{\mathbb D}_{\infty}^{\infty-}.
\end{eqnarray}
By duality, for every $\alpha\geqslant 0$ and $1<p,q,r<\infty$ such that $1/p+1/q=1/r$, there exists a positive
constant $C=C(p,q;\alpha)$ such that
\begin{eqnarray}
\|F\cdot G\|_{-\alpha,r}\leqslant C\|F\|_{\alpha,p}\cdot\|G\|_{-\alpha,q}\quad\text{for all}\quad F,G\in
{\mathbb D}_{\infty}^{\infty-}.
\end{eqnarray}

The natural coupling of $G\in{\mathbb D}_{\alpha}^p$ and $\Phi\in({\mathbb D}_{\alpha}^p)'={\mathbb D}_{-\alpha}^{p/(p-1)}$
(or that of $G\in{\mathbb D}_{\infty}^{\infty-}$ and $\Phi\in{\mathbb D}_{-\infty}^{1+}$) is denoted by $E(G\cdot\Phi)$;
in particular, when $G=\textbf{1}$, the Wiener functional identically equal to one, $E(\textbf{1}\cdot\Phi)$ is
simply denoted by $E(\Phi)$ and is called the generalized expectation of $\Phi\in{\mathbb D}_{-\infty}^{1+}$.
$E(\Phi)$ coincides with the usual expectation $\int_B\Phi(\omega)\mu(d\omega)$ when $\Phi$ is given by an integrable
random variable.

Throughout the paper we shall use the following classical function spaces on $\mR^d$. For a multiindex $\alpha=(\alpha_1,\cdots,\alpha_k)\in\{1,\cdots,d\}^k$ we write $|\alpha|=k$ for its length and $\partial_{\alpha}f=\partial_{x_{\alpha_1}}\cdots\partial_{x_{\alpha_k}}f$ for the corresponding partial derivative of a sufficiently smooth function $f$. For an integer $m\geq0$ we denote by $C^m(\mR^d)$ the space of all functions $f:\mR^d\to\mR$ with continuous partial derivatives up to order $m$, by $C_b^m(\mR^d)$ the subspace of functions $f\in C^m(\mR^d)$ whose partial derivatives up to order $m$ are bounded, equipped with the norm
\ce
\|f\|_{C_b^m}=\sum_{|\alpha|\leq m}\|\partial_{\alpha}f\|_{\infty},
\de
where $\|\cdot\|_{\infty}$ denotes the (essential) supremum norm, and by $C_c^m(\mR^d)$ the subspace of functions $f\in C^m(\mR^d)$ with compact support. We also set
\ce
C^{\infty}(\mR^d)=\bigcap_{m\geq0}C^m(\mR^d),\qquad C_c^{\infty}(\mR^d)=\bigcap_{m\geq0}C_c^m(\mR^d),
\de
and we denote by $C_p^{\infty}(\mR^d)$ the set of all infinitely differentiable functions
$f:\mR^d\rightarrow\mR$ such that $f$ and all of its partial derivatives have polynomial growth.

In the Malliavin calculus, a key role is played by the Malliavin covariance matrix which is defined as follows.
\begin{definition}
Suppose that $F=(F_1,\cdots,F_d)$ is a random vector whose components belong to the space ${\mathbb D}_1^{\infty-}$.
We associate with $F$ the following random symmetric nonnegative definite matrix:
\begin{eqnarray*}
\Sigma_F(\omega)=(\sigma_F^{ij}(\omega))_{1\leqslant i,j\leqslant d}:=(\langle DF_i,DF_j\rangle_H)_{1\leqslant i,j\leqslant d}.
\end{eqnarray*}
The matrix $\Sigma_F(\omega)$ will be called the Malliavin covariance matrix of the random vector $F$. We will say
that such a random vector $F$ is nondegenerate
if the Malliavin covariance matrix $\Sigma_F(\omega)$ is invertible a.s. and
\begin{eqnarray*}
\Gamma_F(\omega)=(\gamma_F^{ij}(\omega))_{1\leqslant i,j\leqslant d}:=(\Sigma_F(\omega))^{-1}\in\bigcap_{p\geqslant 1}
\mathbb{L}^p({\mathbb R}^d\otimes{\mathbb R}^d).
\end{eqnarray*}
\end{definition}

The fractional order Sobolev spaces on Wiener space were introduced by Watanabe \cite{W}, who also studied
their applications to solutions of stochastic differential equations. An equivalent description
of these spaces by the $K$-method was then given in \cite{AMR} and \cite{HR}. We now recall the related results
that will be needed in the sequel.
\begin{definition}
For $0<\alpha<1$ we define
\begin{eqnarray}
{\mathcal E}^p_{\alpha}=({\mathbb L}^p, {\mathbb D}_1^p)_{\alpha,p},
\end{eqnarray}
where $(\cdot,\cdot)$ denotes the real interpolation space as in \cite{T}.
\end{definition}

There are several equivalent norms on ${\mathcal E}_{\alpha}^p$ (cf. \cite{T}). The
one we shall use is given by Peetre's $K$-method:
\begin{eqnarray}
\|f\|_{\alpha,p}:=
\big[ \int_0^1 [\epsilon^{-\alpha} K(\epsilon, f)]^p
\frac{d\epsilon}{\epsilon} \big]^{\frac{1}{p}},
\end{eqnarray}
where
\begin{eqnarray}
K(\epsilon ,f) = \inf \{ \|f_1\|_p + \epsilon \|f_2\|_{1,p},\; f_1 +
f_2 =f,\; f_1\in {\mathbb L}^p,\; f_2\in{\mathbb D}_1^p \}.
\end{eqnarray}

\begin{remark}\label{k}
Let $0<\alpha<1$. For each integer $n\geq1$, let $K_n(F)=K(2^{-n},F)$. Then $F\in{\mathcal E}_{\alpha}^p$ if and only if
$$
\sum_{n=1}^\infty 2^{np\alpha}|K_n(F)|^p<\infty.
$$
\end{remark}

For the relationship between ${\mathbb D}_{\alpha}^p$ and ${\mathcal E}_{\alpha}^p$, the following result is proved in
\cite[Theorem 1.1]{W} directly using specific properties of Wiener functionals.
\begin{proposition}\label{two types}
For every $0<\alpha<1$, $1<p<\infty$ and $\varepsilon>0$, we have
\begin{eqnarray*}
{\mathcal E}_{\alpha-\varepsilon}^p\subset{\mathbb D}_{\alpha}^p\subset{\mathcal E}_{\alpha+\varepsilon}^p.
\end{eqnarray*}
\end{proposition}

The following proposition can be found in \cite[Theorem 2.9]{HRZ}.
\begin{proposition}\label{inverse}
Let $\alpha=k+\sigma$, $k\in{\mathbb N}$, and $0<\sigma\leqslant 1$. Suppose that $F\in{\mathcal E}_{\alpha}^p$ is a nonnegative Wiener functional,
$G\in {\mathbb D}_{k+1}^r$, and
$A$ is a measurable set such that $G1_A=0$ and
\begin{eqnarray*}
E[|F|^{-p}1_{A^c}]<\infty,\quad\text{for every}\quad p>1.
\end{eqnarray*}
Then for any $\beta<\alpha$, we have
\begin{eqnarray}
\frac{1}{F}\cdot G\in{\mathcal E}_{\beta}^{\infty-}:=\bigcap_{1<p<\infty}{\mathcal E}_{\beta}^p.
\end{eqnarray}
\end{proposition}

Consider on $\mR^d$ the Poisson kernel $Q_d$, i.e. the fundamental solution of $\Delta Q_d=\delta_0$,
\ce
Q_d(x)=
\begin{cases}
\max\{x,0\},&d=1;\\
a_2^{-1}\ln |x|,&d=2;\\
-\dfrac{1}{(d-2)a_d}|x|^{2-d},&d>2,
\end{cases}
\de
where $a_d$ denotes the area of the unit sphere of $\mR^d$. Their gradients are given by
\ce
\partial_iQ_d(x)=a_d^{-1}\frac{x_i}{|x|^d},\qquad i=1,\cdots,d,\quad d\geq2,
\de
whereas $\partial_1Q_1=1_{(0,\infty)}$ when $d=1$. In particular, $\partial_iQ_d$ is homogeneous of degree $-(d-1)$ and $\partial_j\partial_iQ_d$ is homogeneous of degree $-d$.
Malliavin and Thalmaier \cite{MT} give the following formula: for any $f\in C_c^1(\mR^d)$,
\ce
f=(\nabla Q_d)*\nabla f.
\de
This representation of $f$ is called the Riesz transform of $f$, and it plays an important role in the representation formula for conditional expectations.

Throughout the paper, $C$ will denote different constants
(which may or may not depend on the indices) whose exact values are not important.

\section{Fractional order Sobolev spaces and density estimates}
\subsection{Fractional order Sobolev spaces associated with a probability measure}
We now introduce the fractional order Sobolev spaces associated with a probability measure.
We consider a probability measure $\mu$ on $\mR^d$ (with the Borel $\sigma$-field). We denote by $L_{\mu}^p$ the space $\{\phi;\int|\phi(x)|^p\mu(dx)<\infty\}$, equipped with the norm $\|\phi\|_{L_{\mu}^p}=(\int|\phi(x)|^p\mu(dx))^{1/p}$.
\bd
We denote by $W_{\mu}^{1,p}$ the space of functions $\phi\in L_{\mu}^p$ for which there exist some functions $\theta_i\in L_{\mu}^p$, $i=1,\cdots,d$ such that, for every test function $f\in C_c^{\infty}(\mR^d)$, one has
\ce
\int\partial_if(x)\phi(x)\mu(dx)=-\int f(x)\theta_i(x)\mu(dx),
\de
with the norm
\ce
\|\phi\|_{W_{\mu}^{1,p}}=\|\phi\|_{L_{\mu}^p}+\sum_{i=1}^d\|\partial_i^{\mu}\phi\|_{L_{\mu}^p},
\de
where $\partial_i^{\mu}\phi=\theta_i$.

Similarly, we define the Sobolev spaces of higher orders. Recall that for a multiindex $\alpha=(\alpha_1,\cdots,\alpha_m)\in\{1,\cdots,d\}^m$ and a function $f\in C^m(\mR^d)$, $\partial_{\alpha}f=\partial_{x_{\alpha_1}}\cdots\partial_{x_{\alpha_m}}f$ denotes the partial derivative corresponding to $\alpha$ (see Section 2). We define $W_{\mu}^{m,p}$ to be the space of functions $\phi\in L_{\mu}^p$ such that for every multiindex $\alpha$ with $|\alpha|\leq m$ there exist some functions $\theta_{\alpha}\in L_{\mu}^p$ such that
\ce
\int\partial_{\alpha}f(x)\phi(x)\mu(dx)=(-1)^{|\alpha|}\int f(x)\theta_{\alpha}(x)\mu(dx)\quad\forall f\in C_c^{\infty}(\mR^d)
\de
with the norm
\ce
\|\phi\|_{W_{\mu}^{m,p}}=\|\phi\|_{L_{\mu}^p}+\sum_{|\alpha|\leq m}\|\partial_{\alpha}^{\mu}\phi\|_{L_{\mu}^p},
\de
where $\partial_{\alpha}^{\mu}\phi=\theta_{\alpha}$.

We define $W_{\mu}^{-k,p}=(W_{\mu}^{k,q})'$ where $1/p+1/q=1$. Let $1<p<\infty$, $\alpha\in{\mathbb R}$. We also define fractional order Sobolev spaces associated with a probability measure on ${\mathbb R}^d$ as follows:
\begin{eqnarray*}
{\mathcal E}_{\mu}^{\alpha,p}=(W_{\mu}^{m-1,p},W_{\mu}^{m,p})_{\theta,p}
\quad\text{where $m$ is the smallest integer larger than $\alpha$}
\end{eqnarray*}
equipped with the norm
\begin{eqnarray*}
\|u\|_{{\mathcal E}_{\mu}^{\alpha,p}}=[\int_0^{1}|\epsilon^{-\theta}K(\epsilon,u)|^p\frac{d\epsilon}{\epsilon}]^{1/p},
\end{eqnarray*}
where
\begin{eqnarray*}
K(\epsilon ,f)=\inf\{\|f_1\|_{W_{\mu}^{m-1,p}}+\epsilon \|f_2\|_{W_{\mu}^{m,p}},\; f_1+f_2=f,\; f_1\in W_{\mu}^{m-1,p},\; f_2\in W_{\mu}^{m,p}\}
\end{eqnarray*}
and $\theta:=\alpha-(m-1)$.

We use the notation $L^p$, $W^{m,p}$ and $\cE^{\alpha,p}$ for the corresponding spaces associated with the Lebesgue measure, i.e. the standard $L^p$, Sobolev and fractional-order Sobolev spaces used in the literature. For $0<\delta\leqslant 1$ and $1<p<\infty$ we shall simply write
\ce
W_{\mu}^{\delta,p}:=\cE_{\mu}^{\delta,p}=(L_{\mu}^p,W_{\mu}^{1,p})_{\delta,p}.
\de
If $D\subset\mR^d$ is an open set, we denote by $W_{\mu}^{m,p}(D)$ the space of functions $\phi$ which satisfy the integration by parts formula for test functions $f$ with compact support contained in $D$. The analogous localized spaces $W^{m,p}(D)$, $\cE_{\mu}^{\alpha,p}(D)$, $L^p(D)$ and $L_{\mu}^p(D)$ are defined in the same way.
\ed

To treat fractional orders we also need the following Besov-type spaces associated with $\mu$, which are built on the \emph{ordinary} weak derivatives instead of the integration by parts structure.
\bd\label{plainbesov}
Let $1<q<\infty$. We denote by $X_{\mu}^{1,q}$ the space of functions $u\in L_{\mu}^q$ such that the ordinary weak derivatives $\partial_iu$, $i=1,\cdots,d$, belong to $L_{\mu}^q$, equipped with the norm
\ce
\|u\|_{X_{\mu}^{1,q}}=\|u\|_{L_{\mu}^q}+\sum_{i=1}^d\|\partial_iu\|_{L_{\mu}^q}.
\de
For $0<s<1$ we define
\ce
B_{\mu}^{s,q}:=(L_{\mu}^q,X_{\mu}^{1,q})_{s,q}.
\de
\ed
For $0<s<1$ and $1<q<\infty$ we denote by $B^{s,q}(\mR^d)$ the classical fractional Sobolev (Slobodeckij) space associated with the Lebesgue measure, equipped with the norm
\ce
\|u\|_{B^{s,q}}^q=\|u\|_{L^q}^q+\int_{\mR^d}\int_{\mR^d}\frac{|u(y)-u(z)|^q}{|y-z|^{d+sq}}\,dy\,dz.
\de
More generally, for an open set $D\subset\mR^d$ we denote by $B^{s,q}(D)$ the corresponding space on $D$, equipped with the norm
\ce
\|u\|_{B^{s,q}(D)}^q=\|u\|_{L^q(D)}^q+[u]_{B^{s,q}(D)}^q,\qquad
[u]_{B^{s,q}(D)}^q=\int_D\int_D\frac{|u(y)-u(z)|^q}{|y-z|^{d+sq}}\,dy\,dz,
\de
with the convention $B^{0,q}(D)=L^q(D)$.
\br\label{bbound}
If $\mu(dx)=p_{\mu}(x)dx$ with $\|p_{\mu}\|_{\infty}<\infty$, then for every $u\in B^{s,q}(\mR^d)$,
\ce
\|u\|_{B_{\mu}^{s,q}}\leqslant \|p_{\mu}\|_{\infty}^{1/q}\|u\|_{B^{s,q}},
\de
which follows by interpolating the obvious inequalities $\|u\|_{L_{\mu}^q}\leqslant\|p_{\mu}\|_{\infty}^{1/q}\|u\|_{L^q}$ and $\|u\|_{X_{\mu}^{1,q}}\leqslant\|p_{\mu}\|_{\infty}^{1/q}\|u\|_{W^{1,q}}$.
\er

For a function $\phi\in L_{\mu}^1$, we denote by $\mu_{\phi}$ the signed finite measure defined by
\ce
\mu_{\phi}(dx):=\phi(x)\mu(dx).
\de

In this subsection we state the main abstract results concerning the density of a probability measure satisfying the fractional Sobolev condition $1\in W_{\mu}^{\delta,p}$ with $\delta p>d$. The following result gives a representation formula for the density and its derivatives.
\bt\label{elementary}
Let $1<p<\infty$, $0<\delta\leqslant 1$ with $\delta p>d$, and let $1\in W_{\mu}^{\delta,p}$, $\phi\in W_{\mu}^{\delta,p}$. Then $\mu_{\phi}(dx)=p_{\mu_{\phi}}(x)dx$ with
\be\label{representation1}
p_{\mu_{\phi}}(x)=-\sum_{i=1}^d\langle\partial_i^{\mu}\phi,\partial_iQ_d(\cdot-x)\rangle,
\ee
where the pairing is the one to be defined in Lemma \ref{pairing}, and $p_{\mu_{\phi}}$ is continuous on $\mR^d$. Moreover, if $m\geq1$ and $\phi\in W_{\mu}^{m+\delta,p}$, then
\be\label{representation2}
\partial_{\alpha}p_{\mu_{\phi}}(x)=-\sum_{i=1}^d\langle\partial_{(\alpha,i)}^{\mu}\phi,\partial_iQ_d(\cdot-x)\rangle,
\ee
where $\alpha$ is any multiindex of length less than or equal to $m$, and the following alternative representation formula holds:
\be\label{alternative}
\partial_{\alpha}p_{\mu_{\phi}}(x)=p_{\mu}(x)\partial_{\alpha}^{\mu}\phi(x).
\ee
\et

Now we discuss some regularity properties of the density.
\bt\label{regularity}
Let $1<p<\infty$, $0<\delta'<\delta\leqslant 1$ with $\delta p>d$, $m\geq 0$ and let $1\in W_{\mu}^{\delta,p}$, $\phi\in W_{\mu}^{m+\delta,p}$. Then
\begin{enumerate}
  \item We have $p_{\mu_{\phi}}\in W^{m+\delta',p}$ and
  \ce
  \|p_{\mu_{\phi}}\|_{W^{m+\delta',p}}\leq (C(\delta,p,d,\|1\|_{W_{\mu}^{\delta,p}})\|1\|_{W_{\mu}^{\delta,p}})^{(p-1)/p}\|\phi\|_{W_{\mu}^{m+\delta,p}}.
  \de
  Moreover, for any multiindex $\alpha$ such that $0\leq|\alpha|=l\leq m$, we have
  \ce
  \|\partial_{\alpha}p_{\mu_{\phi}}\|_{\infty}\leq C(\delta,p,d,\|1\|_{W_{\mu}^{\delta,p}})\|1\|_{W_{\mu}^{\delta,p}}\|\phi\|_{W_{\mu}^{m+\delta,p}}.
  \de
  \item Assume in addition that $m+\delta'>d/p$. Let $k_1=[m+\delta'-d/p]$ be the integer part of $m+\delta'-d/p$ and $k_2=\{m+\delta'-d/p\}$ the fractional part. If $k_2>0$, we have $p_{\mu_{\phi}}\in C^{k_1}$; moreover, for any multiindex $\beta$ such that $0\leq|\beta|=k\leq k_1-1$, $\partial_{\beta}p_{\mu_{\phi}}$ is Lipschitz continuous, and for any multiindex $\beta$ such that $|\beta|=k_1$, $\partial_{\beta}p_{\mu_{\phi}}$ is H\"{o}lder continuous of exponent $k_2$. In the borderline case $k_2=0$ (when $m+\delta'-d/p$ is an integer), the usual Sobolev embedding gives $p_{\mu_{\phi}}\in C^{k_1-1,\alpha}$ for every $\alpha<1$ instead. This exceptional case is immaterial for the applications below, since $\delta'<\delta$ and $p$ can always be chosen so that $m+\delta'-d/p$ is not an integer.
\end{enumerate}
\et

Taking $\phi=1$ in Theorems \ref{elementary} and \ref{regularity}, we obtain the following result concerning the density of $\mu$ itself.
\bc\label{holderdensity}
Let $1<p<\infty$, $0<\delta\leqslant 1$ with $\delta p>d$, and let $1\in W_{\mu}^{\delta,p}$. Then $\mu$ is absolutely continuous with respect to the Lebesgue measure, $\mu(dx)=p_{\mu}(x)dx$, and the density $p_{\mu}$ admits the representation
\ce
p_{\mu}(x)=-\sum_{i=1}^d\langle\partial_i^{\mu}1,\partial_iQ_d(\cdot-x)\rangle,
\de
where the pairing is the one defined in Lemma \ref{pairing}. Moreover, $p_{\mu}$ is H\"{o}lder continuous of every exponent $\beta<\delta-d/p$, and there exists a positive constant $C(\delta,p,d,\|1\|_{W_{\mu}^{\delta,p}})$ such that
\ce
\|p_{\mu}\|_{\infty}\leq C(\delta,p,d,\|1\|_{W_{\mu}^{\delta,p}})\|1\|_{W_{\mu}^{\delta,p}}.
\de
\ec

We now establish the auxiliary lemmas needed in the proofs of the above results.
The following lemma gives a rigorous meaning to the ``fractional derivative'' $\partial_i^{\mu}\phi$ acting on singular kernels, through a bilinear pairing.
\bl\label{pairing}
Let $1<p<\infty$, $1/p+1/q=1$ and $0<\delta\leqslant 1$. Then for every $\phi\in W_{\mu}^{1,p}$ and $k\in C_c^1(\mR^d)$,
\be\label{skew}
\int k\,\partial_i^{\mu}\phi\,\mu(dx)=-\int\phi\,\partial_ik\,\mu(dx).
\ee
Consequently, the bilinear form $B_i(\phi,k):=\int k\,\partial_i^{\mu}\phi\,\mu(dx)$ extends uniquely to a bounded bilinear form on $W_{\mu}^{\delta,p}\times B_{\mu}^{1-\delta,q}$. We denote the extension by
\ce
\langle\partial_i^{\mu}\phi,k\rangle:=B_i(\phi,k),\qquad \phi\in W_{\mu}^{\delta,p},\ k\in B_{\mu}^{1-\delta,q}.
\de
\el
\begin{proof}
For $\phi\in W_{\mu}^{1,p}$ and $k\in C_c^1(\mR^d)$ we have $k\phi\in W_{\mu}^{1,1}$ and $\partial_i^{\mu}(k\phi)=k\partial_i^{\mu}\phi+\phi\partial_ik$; this follows from the definition of $\partial_i^{\mu}$ applied to the test functions $fk$ with $f\in C_c^{\infty}$. Taking test functions $f_n\in C_c^{\infty}$ with $0\leqslant f_n\uparrow 1$ and $|\nabla f_n|\leqslant 2^{-n}$ in the integration by parts formula for $k\phi$ and letting $n\to\infty$ gives $\int\partial_i^{\mu}(k\phi)d\mu=0$, which is exactly (\ref{skew}). Hence
\ce
|B_i(\phi,k)|\leqslant\|\phi\|_{W_{\mu}^{1,p}}\|k\|_{L_{\mu}^q},\qquad
|B_i(\phi,k)|\leqslant\|\phi\|_{L_{\mu}^p}\|k\|_{X_{\mu}^{1,q}},
\de
where the second inequality uses (\ref{skew}). Since $C_c^1(\mR^d)$ is dense in $X_{\mu}^{1,q}$, both estimates hold for all $k\in X_{\mu}^{1,q}$, and bilinear real interpolation between $(W_{\mu}^{1,p},L_{\mu}^q)$ and $(L_{\mu}^p,X_{\mu}^{1,q})$ yields the boundedness of $B_i$ on $(W_{\mu}^{1,p},L_{\mu}^p)_{1-\delta,p}\times(L_{\mu}^q,X_{\mu}^{1,q})_{1-\delta,q}=W_{\mu}^{\delta,p}\times B_{\mu}^{1-\delta,q}$.
\end{proof}

\br\label{dualnorm}
In particular, for $\phi\in W_{\mu}^{\delta,p}$ the element $\partial_i^{\mu}\phi\in(B_{\mu}^{1-\delta,q})^*$ is well defined and
\ce
\|\partial_i^{\mu}\phi\|_{(B_{\mu}^{1-\delta,q})^*}\leqslant C\|\phi\|_{W_{\mu}^{\delta,p}}.
\de
Iterating this estimate, if $\phi\in W_{\mu}^{m+\delta,p}$ for some integer $m\geq0$, then
\ce
\|\partial_{(\alpha,i)}^{\mu}\phi\|_{(B_{\mu}^{1-\delta,q})^*}\leqslant C\|\phi\|_{W_{\mu}^{m+\delta,p}}
\de
for every multiindex $\alpha$ with $|\alpha|\leq m$ and every $i=1,\cdots,d$.
\er

We now estimate the Riesz kernels in the above Besov-type spaces. This is the key point where the condition $\delta p>d$ comes in. We fix once and for all a cutoff function $\chi\in C_c^{\infty}(\mR^d)$ such that $0\leq\chi\leq1$, $\chi=0$ on $B(0,1)$ and $\chi=1$ on $\mR^d\setminus B(0,2)$; for $x\in\mR^d$ and $\rho>0$ we write $\chi_{x,\rho}(y):=\chi((y-x)/\rho)$.
\bl\label{kernel}
Let $1<p<\infty$, $0<\delta\leqslant 1$ with $\delta p>d$, and set $q=p/(p-1)$. Then there exist constants $C(\delta,p,d)$ and $C(p,d)$ such that, with
\ce
M_0(\rho):=C(\delta,p,d)\rho^{(\delta p-d)/p},\qquad M_1(\rho):=C(p,d)(\rho^{-(d-1)}+\rho^{-d}),
\de
for every $x\in\mR^d$, $0<\rho\leqslant1$ and $i=1,\cdots,d$:
\begin{enumerate}
\item[(i)] $\|\partial_iQ_d(\cdot-x)\|_{B^{1-\delta,q}(B(x,\rho))}\leqslant M_0(\rho)$;
\item[(ii)] $\|\chi_{x,\rho}\,\partial_iQ_d(\cdot-x)\|_{X_{\mu}^{1,q}}\leqslant M_1(\rho)$;
\item[(iii)] if, in addition, $\mu(dx)=p_{\mu}(x)dx$ with $\|p_{\mu}\|_{\infty}<\infty$, then
\ce
\Gamma_p(\mu):=\sup_{x\in\mR^d}\sum_{i=1}^d\|\partial_iQ_d(\cdot-x)\|_{B_{\mu}^{1-\delta,q}}\leqslant M_1(\rho)+\|p_{\mu}\|_{\infty}^{1/q}M_0(\rho).
\de
\end{enumerate}
\el
\begin{proof}
We shall use the norms $\|\cdot\|_{B^{1-\delta,q}(D)}$ on the domains $D=B(x,\rho)$ and $D=\mR^d$ defined after Definition \ref{plainbesov} (the convention $B^{0,q}=L^q$ covers the case $\delta=1$).

(i) By translation invariance it suffices to take $x=0$. For $d\geq2$, $\partial_iQ_d(z)=a_d^{-1}z_i|z|^{-d}$ is smooth away from the origin and homogeneous of degree $-(d-1)$; for $d=1$, $\partial_1Q_1=1_{(0,\infty)}$. The dilation $y=\rho\,\widetilde y$ gives
\ce
\|\partial_iQ_d\|_{L^q(B(0,\rho))}^q=\rho^{\,d-(d-1)q}\|\partial_iQ_d\|_{L^q(B(0,1))}^q,\\
{}[\partial_iQ_d]_{B^{1-\delta,q}(B(0,\rho))}^q=\rho^{\,d-(d-\delta)q}[\partial_iQ_d]_{B^{1-\delta,q}(B(0,1))}^q.
\de
The first integral is finite if and only if $(d-1)q<d$, which is equivalent to $p>d$ and follows from $\delta p>d$ since $\delta\leqslant1$. For the seminorm, the standard integrability criterion for radial singularities gives
\ce
[\partial_iQ_d]_{B^{1-\delta,q}(B(0,1))}<\infty\quad\Longleftrightarrow\quad (d-\delta)q<d.
\de
Indeed, this can be obtained by splitting the double integral according as $|y-z|\geq|y|/2$ or $|y-z|<|y|/2$: on the first region one uses
$|\partial_iQ_d(y)-\partial_iQ_d(z)|^q\leq C(|y|^{-(d-1)q}+|z|^{-(d-1)q})$ and integrates by polar coordinates, while on the second region the mean value theorem together with $|\partial_j\partial_iQ_d(w)|\leqslant C|w|^{-d}$ gives
$|\partial_iQ_d(y)-\partial_iQ_d(z)|\leqslant C|y-z|\,|y|^{-d}$; both contributions are bounded by
\ce
C\int_0^1 r^{\,d-1-(d-\delta)q}\,dr,
\de
which is finite exactly when $(d-\delta)q<d$, equivalently $\delta p>d$. For $d=1$, $\partial_1Q_1$ is a step function with a single jump at $0$, and the same computation gives $(1-\delta)q<1$, i.e. $\delta p>1$.

Since
\ce
d-(d-1)q=\frac{p-d}{p-1},\qquad d-(d-\delta)q=\frac{\delta p-d}{p-1},
\de
taking $q$-th roots yields
\ce
\|\partial_iQ_d(\cdot-x)\|_{L^q(B(x,\rho))}\leqslant C\rho^{(p-d)/p},\qquad
[\partial_iQ_d(\cdot-x)]_{B^{1-\delta,q}(B(x,\rho))}\leqslant C\rho^{(\delta p-d)/p}.
\de
As $\delta\leqslant1$ implies $(p-d)/p\geq(\delta p-d)/p$, and $0<\rho\leqslant1$ implies $\rho^{(p-d)/p}\leqslant\rho^{(\delta p-d)/p}$, the $L^q$ contribution is dominated by the Besov seminorm contribution, and the whole norm on $B(x,\rho)$ is bounded by $M_0(\rho)$. This proves (i).

(ii) Define the exterior kernel
\ce
h_{x,\rho}(y):=\chi_{x,\rho}(y)\,\partial_iQ_d(y-x),
\de
which vanishes on $B(x,\rho)$. On $\{|y-x|\geq\rho\}$ one has
\ce
|\partial_iQ_d(y-x)|\leqslant C\rho^{-(d-1)},\qquad |\partial_j\partial_iQ_d(y-x)|\leqslant C\rho^{-d}.
\de
Since
\ce
\partial_jh_{x,\rho}=\chi_{x,\rho}\,\partial_j\partial_iQ_d(\cdot-x)+(\partial_j\chi_{x,\rho})\,\partial_iQ_d(\cdot-x)
\de
and $|\partial_j\chi_{x,\rho}|\leqslant C\rho^{-1}$, we obtain, using $\mu(\mR^d)=1$,
\ce
\|h_{x,\rho}\|_{L_{\mu}^q}\leqslant\|h_{x,\rho}\|_{\infty}\leqslant C\rho^{-(d-1)},\qquad
\|\partial_jh_{x,\rho}\|_{L_{\mu}^q}\leqslant\|\partial_jh_{x,\rho}\|_{\infty}\leqslant C\rho^{-d}.
\de
Summing over $j$ gives $\|h_{x,\rho}\|_{X_{\mu}^{1,q}}\leqslant C(p,d)(\rho^{-(d-1)}+\rho^{-d})=M_1(\rho)$.

(iii) Define the interior kernel
\ce
g_{x,\rho}(y):=(1-\chi_{x,\rho}(y))\,\partial_iQ_d(y-x),
\de
which is supported in $B(x,2\rho)$. We claim that
\be\label{interior-cutoff}
\|g_{x,\rho}\|_{B^{1-\delta,q}(\mR^d)}\leqslant C(\delta,p,d)\rho^{(\delta p-d)/p}=M_0(\rho),\qquad 0<\rho\leqslant1.
\ee
Indeed, $g_{x,\rho}(y)=\rho^{-(d-1)}g_{0,1}((y-x)/\rho)$, where $g_{0,1}$ is a fixed compactly supported function which equals $\partial_iQ_d$ in a neighborhood of $0$; by the criterion proved in (i), $g_{0,1}\in B^{1-\delta,q}(\mR^d)$ exactly under $\delta p>d$. The dilation formulas for the Besov norm then give (\ref{interior-cutoff}), with the $L^q$ contribution dominated by the seminorm contribution because $\rho\leqslant1$.

The cutoff now gives the exact decomposition $\partial_iQ_d(\cdot-x)=g_{x,\rho}+h_{x,\rho}$, with no boundary term. By Remark \ref{bbound} and (\ref{interior-cutoff}),
\ce
\|g_{x,\rho}\|_{B_{\mu}^{1-\delta,q}}\leqslant\|p_{\mu}\|_{\infty}^{1/q}\|g_{x,\rho}\|_{B^{1-\delta,q}}\leqslant\|p_{\mu}\|_{\infty}^{1/q}M_0(\rho).
\de
On the other hand, since $B_{\mu}^{1-\delta,q}=(L_{\mu}^q,X_{\mu}^{1,q})_{1-\delta,q}$, the definition of the interpolation norm via the $K$-functional gives
\ce
\|u\|_{B_{\mu}^{1-\delta,q}}\leqslant C(\delta,p)\big(\|u\|_{L_{\mu}^q}+\|u\|_{X_{\mu}^{1,q}}\big),
\de
and hence $\|h_{x,\rho}\|_{B_{\mu}^{1-\delta,q}}\leqslant C M_1(\rho)$. Adding the two estimates, summing over $i=1,\cdots,d$ and taking the supremum over $x$ proves (iii).
\end{proof}

The result stated below is borrowed from \cite[Lemma 3]{BC1}.
\bl\label{t1}
Let $p_n$, $n\in\mN$, be a sequence of probability densities such that $\sup_n\|p_n\|_{\infty}=C_{\infty}<\infty$. Suppose also that the sequence of probability measures $\mu_n(dx)=p_n(x)dx$, $n\in\mN$ converges weakly to a probability measure $\mu$. Then $\mu(dx)=p(x)dx$ and $\|p\|_{\infty}\leq C_{\infty}$.
\el

If $\psi$ is an integrable function on $\mR^d$ and $\mu$ is a probability measure on $\mR^d$, we denote by $\psi*\mu$ the convolution of $\psi$ with $\mu$, that is, the finite signed measure on $\mR^d$ defined by
\ce
\int_{\mR^d} f(x)(\psi*\mu)(dx)=\int_{\mR^d}\int_{\mR^d} f(x+y)\psi(y)\,dy\,\mu(dx)
\de
for every bounded measurable function $f$. When $\psi\geq0$ and $\int\psi=1$, $\psi*\mu$ is the law of $F+\xi$, where $F\sim\mu$ and $\xi\sim\psi(x)dx$ are independent random variables; it is then absolutely continuous with respect to the Lebesgue measure, with density $x\mapsto\int_{\mR^d}\psi(x-y)\mu(dy)$.

The next lemma shows that conditioning on the independent convolution variable does not increase the fractional Sobolev norm.
\bl\label{t2}
Let $p\geq 1$ and $0\leq\alpha\leq1$. Let $F\sim\mu$ and $\xi\sim\psi(x)dx$ be independent random variables, with $\psi\geq0$ and $\int\psi=1$, and set $T\phi(x):=E[\phi(F)|F+\xi=x]$. Then $T$ maps $W_{\mu}^{\alpha,p}$ into $W_{\psi*\mu}^{\alpha,p}$ and
\ce
\|T\phi\|_{W_{\psi*\mu}^{\alpha,p}}\leq\|\phi\|_{W_{\mu}^{\alpha,p}}.
\de
In particular $T1=1$, so if $1\in W_{\mu}^{\alpha,p}$ then $1\in W_{\psi*\mu}^{\alpha,p}$ with $\|1\|_{W_{\psi*\mu}^{\alpha,p}}\leq\|1\|_{W_{\mu}^{\alpha,p}}$.
\el
\begin{proof}
The case $\alpha=0$ is trivial and the case $\alpha=1$ is exactly \cite[Lemma 3]{BC1}. We therefore assume $0<\alpha<1$ for simplicity. We first prove that
\be\label{tk}
K_{\psi*\mu}(\epsilon,T\phi)\leq K_{\mu}(\epsilon,\phi),\quad\forall\epsilon>0,
\ee
where $K_{\mu}$ denotes the $K$-functional associated with the pair $(L_{\mu}^p,W_{\mu}^{1,p})$ and $K_{\psi*\mu}$ the one associated with the pair $(L_{\psi*\mu}^p,W_{\psi*\mu}^{1,p})$.

Let $f_1+f_2=\phi$ with $f_1\in L_{\mu}^p$ and $f_2\in W_{\mu}^{1,p}$ be arbitrary, and set
\ce
g_i(x)=E[f_i(F)|F+\xi=x],\quad i=1,2.
\de
Then $g_1+g_2=T\phi$ and by Jensen's inequality for conditional expectations,
\ce
\|g_1\|_{L_{\psi*\mu}^p}^p=E[|E(f_1(F)|F+\xi)|^p]\leq E[|f_1(F)|^p]=\|f_1\|_{L_{\mu}^p}^p.
\de
Moreover, for every $f\in C_c^{\infty}(\mR^d)$ one has
\ce
-\int\partial_if(x)g_2(x)(\psi*\mu)(dx)&=&-E[\partial_if(F+\xi)g_2(F+\xi)]\\
&=&-E[\partial_if(F+\xi)f_2(F)]\\
&=&-E[f_2(F)E(\partial_if(F+\xi)|F)]\\
&=&-E[f_2(F)\int\psi(y)\partial_if(F+y)dy]\\
&=&E[\int\psi(y)f(F+y)dy\cdot\partial_i^{\mu}f_2(F)]\\
&=&E[f(F+\xi)\partial_i^{\mu}f_2(F)]\\
&=&\int f(x)E[\partial_i^{\mu}f_2(F)|F+\xi=x](\psi*\mu)(dx),
\de
where, in the fifth line, we have used the definition of $\partial_i^{\mu}f_2$ with the test function $x\mapsto f(x+y)$ (which has compact support) and Fubini's theorem. Hence $g_2\in W_{\psi*\mu}^{1,p}$ and
\ce
\partial_i^{\psi*\mu}g_2(x)=E[\partial_i^{\mu}f_2(F)|F+\xi=x].
\de
By Jensen's inequality again, we have
\ce
\|g_2\|_{W_{\psi*\mu}^{1,p}}\leq\|f_2\|_{W_{\mu}^{1,p}}.
\de
Therefore
\ce
K_{\psi*\mu}(\epsilon,T\phi)\leq\|g_1\|_{L_{\psi*\mu}^p}+\epsilon\|g_2\|_{W_{\psi*\mu}^{1,p}}\leq\|f_1\|_{L_{\mu}^p}+\epsilon\|f_2\|_{W_{\mu}^{1,p}},
\de
and (\ref{tk}) follows by taking the infimum over all decompositions $f_1+f_2=\phi$ on the right-hand side. Finally, by the definition of the fractional order Sobolev spaces,
\ce
\|T\phi\|_{W_{\psi*\mu}^{\alpha,p}}&=&\big[\int_0^1|\epsilon^{-\alpha}K_{\psi*\mu}(\epsilon,T\phi)|^p\frac{d\epsilon}{\epsilon}\big]^{1/p}\\
&\leq&\big[\int_0^1|\epsilon^{-\alpha}K_{\mu}(\epsilon,\phi)|^p\frac{d\epsilon}{\epsilon}\big]^{1/p}\\
&=&\|\phi\|_{W_{\mu}^{\alpha,p}},
\de
which completes the proof.
\end{proof}

The key boundedness estimate required for the proofs of Theorem \ref{elementary}, Theorem \ref{regularity}, and Corollary \ref{holderdensity} is provided by the following lemma, which also gives a uniform estimate for the Riesz kernels.
\bl\label{boundmain}
Let $1<p<\infty$, $0<\delta\leqslant 1$ with $\delta p>d$, and let $\mu$ be a probability measure on $\mR^d$ such that $1\in W_{\mu}^{\delta,p}$. Then $\mu(dx)=p_{\mu}(x)dx$ and we have
\be\label{i1}
\|p_{\mu}\|_{\infty}\leq C_1(\delta,p,d,\|1\|_{W_{\mu}^{\delta,p}})\|1\|_{W_{\mu}^{\delta,p}}
\ee
and, with $\Gamma_p(\mu)$ as defined in Lemma \ref{kernel} (iii),
\be\label{i2}
\Gamma_p(\mu)\leq C_2(\delta,p,d,\|1\|_{W_{\mu}^{\delta,p}})\|1\|_{W_{\mu}^{\delta,p}}.
\ee
\el

\begin{proof}
	Set $q=p/(p-1)$. We first prove the lemma under the supplementary assumption that $\mu(dx)=p_{\mu}(x)dx$ with $\|p_{\mu}\|_{\infty}<\infty$. By Step~1 of the proof of Theorem \ref{elementary} below (which, under this supplementary assumption, relies only on Lemma \ref{pairing} and Lemma \ref{kernel} and does not use the present lemma), Lemma \ref{pairing} (cf. Remark \ref{dualnorm}) and Lemma \ref{kernel} (iii), for every $x\in\mR^d$ and every $0<\rho\leq1$,
	\ce
	p_{\mu}(x)&\leqslant&\sum_{i=1}^d\big|\langle\partial_i^{\mu}1,\partial_iQ_d(\cdot-x)\rangle\big|\\
	&\leqslant& C\|1\|_{W_{\mu}^{\delta,p}}\sum_{i=1}^d\|\partial_iQ_d(\cdot-x)\|_{B_{\mu}^{1-\delta,q}}\\
	&\leqslant& C\|1\|_{W_{\mu}^{\delta,p}}\big(M_1(\rho)+\|p_{\mu}\|_{\infty}^{1/q}M_0(\rho)\big).
	\de
	Since $M_0(\rho)=C(\delta,p,d)\rho^{(\delta p-d)/p}\to0$ as $\rho\downarrow0$ (here the condition $\delta p>d$ is used), we can choose $\rho^*\in(0,1]$ such that
	\ce
	C\|1\|_{W_{\mu}^{\delta,p}}M_0(\rho^*)\leq\frac{1}{2};
	\de
	namely, $\rho^*=1$ if the left-hand side already does not exceed $\frac12$ at $\rho=1$, and otherwise $\rho^*$ is the unique solution in $(0,1)$ of the corresponding equality.
	Taking the supremum over $x\in\mR^d$ and using Young's inequality with exponent $q>1$, we obtain
	\ce
	\|p_{\mu}\|_{\infty}\leqslant C\|1\|_{W_{\mu}^{\delta,p}}M_1(\rho^*)+\frac12\|p_{\mu}\|_{\infty}^{1/q}\leqslant C\|1\|_{W_{\mu}^{\delta,p}}M_1(\rho^*)+\frac12\|p_{\mu}\|_{\infty}+C,
	\de
	and hence
	\be\label{eq:pbound}
	\|p_{\mu}\|_{\infty}\leqslant 2C\|1\|_{W_{\mu}^{\delta,p}}M_1(\rho^*)+2C\leqslant C_1(\delta,p,d,\|1\|_{W_{\mu}^{\delta,p}})\|1\|_{W_{\mu}^{\delta,p}}.
	\ee
	Substituting (\ref{eq:pbound}) back into Lemma \ref{kernel} (iii) with $\rho=\rho^*$, we also have
	\ce
	\Gamma_p(\mu)\leqslant M_1(\rho^*)+\|p_{\mu}\|_{\infty}^{1/q}M_0(\rho^*)\leqslant C_2(\delta,p,d,\|1\|_{W_{\mu}^{\delta,p}})\|1\|_{W_{\mu}^{\delta,p}}.
	\de
	Thus the lemma is proved under the supplementary assumption $\|p_{\mu}\|_{\infty}<\infty$.
	
	We now remove this assumption as in the proof of \cite[Theorem 5]{BC1}. Let $\psi\geqslant0$ be a continuous function with $\int\psi=1$ and $\psi(x)=0$ for $|x|\geqslant1$, set $\psi_n(x)=n^d\psi(nx)$ and $\mu_n=\psi_n*\mu$. Then $\mu_n(dx)=p_n(x)dx$ with $p_n$ bounded and continuous. By Lemma \ref{t2}, $1\in W_{\mu_n}^{\delta,p}$ with $\|1\|_{W_{\mu_n}^{\delta,p}}\leqslant\|1\|_{W_{\mu}^{\delta,p}}<\infty$, so the first part of the proof gives
	\ce
	\|p_n\|_{\infty}\leqslant C_1(\delta,p,d,\|1\|_{W_{\mu}^{\delta,p}})\|1\|_{W_{\mu}^{\delta,p}},
	\de
	uniformly in $n$. Since $\mu_n\to\mu$ weakly, Lemma \ref{t1} yields $\mu(dx)=p_{\mu}(x)dx$ with $\|p_{\mu}\|_{\infty}$ bounded by the same constant. Applying the first part of the proof to $\mu$ itself then completes the proof.
\end{proof}

We are now in a position to prove the results stated above.
\begin{proof}[Proof of Theorem \ref{elementary}]
\textbf{Step 1.} We first prove (\ref{representation1}) under the additional assumption that $\mu(dx)=p_{\mu}(x)dx$ with $\|p_{\mu}\|_{\infty}<\infty$. Let $f\in C_c^1(\mR^d)$ with $\mathrm{supp}\,f\subset B_R(0)$ for some $R>1$. Since $f=\Delta(Q_d*f)=\sum_{i=1}^d(\partial_iQ_d)*(\partial_if)$, we have
\ce
\int f d\mu_{\phi}=\sum_{i=1}^d\int\mu(dx)\phi(x)\int\partial_iQ_d(z)\partial_if(x-z)dz.
\de
Therefore
\ce
\int f d\mu_{\phi}&=&\sum_{i=1}^d\int\partial_iQ_d(z)\Big(\int\phi(x)\partial_if(x-z)\mu(dx)\Big)dz\\
&=&-\sum_{i=1}^d\int\partial_iQ_d(z)B_i(\phi,f(\cdot-z))dz\\
&=&-\sum_{i=1}^d B_i\Big(\phi,\int\partial_iQ_d(z)f(\cdot-z)dz\Big)\\
&=&-\sum_{i=1}^d\int f(y)B_i(\phi,\partial_iQ_d(\cdot-y))dy,
\de
which gives (\ref{representation1}). The first equality follows from Fubini's theorem, which applies here since, as in the proof of \cite[Theorem 1]{BC1}, setting $C_R(x):=\{z:(|x|-R)\vee0\leqslant|z|\leqslant|x|+R\}$ one has $\int_{C_R(x)}|\partial_iQ_d(z)|dz\leqslant C_dR$, hence
\ce
\int\int|\phi(x)\partial_iQ_d(z)\partial_if(x-z)|dz\mu(dx)\leqslant C_dR\|\partial_if\|_{\infty}\|\phi\|_{L_{\mu}^1}<\infty;
\de
in the second equality we used (\ref{skew}); the third equality is justified since, by Lemma \ref{pairing} and Remark \ref{bbound},
\ce
\int|\partial_iQ_d(z)|\,\|f(\cdot-z)\|_{B_{\mu}^{1-\delta,q}}dz\leqslant C\|f\|_{C_b^1}(1+\|p_{\mu}\|_{\infty}^{1/q})\int_{C_R(0)}|\partial_iQ_d(z)|dz<\infty;
\de
and the last one since $\sup_{y}\|\partial_iQ_d(\cdot-y)\|_{B_{\mu}^{1-\delta,q}}\leqslant\Gamma_p(\mu)<\infty$ by Lemma \ref{kernel} (iii). The continuity of $x\mapsto p_{\mu_{\phi}}(x)$ follows from the cutoff decomposition used in Lemma \ref{kernel}: the exterior parts $\chi_{x,\rho}\,\partial_iQ_d(\cdot-x)$ depend continuously on $x$ in $X_{\mu}^{1,q}$ by dominated convergence, while the interior parts are bounded in $B_{\mu}^{1-\delta,q}$ uniformly in $x$ by $C\|p_{\mu}\|_{\infty}^{1/q}\rho^{(\delta p-d)/p}$, which tends to zero as $\rho\downarrow0$.

\textbf{Step 2.} The general case, i.e. without the a priori assumption $\|p_{\mu}\|_{\infty}<\infty$, follows from Step 1 and Lemma \ref{boundmain}: the proof of Lemma \ref{boundmain} only uses Step 1 applied to the mollified measures $\mu_n=\psi_n*\mu$, which have bounded continuous densities, and shows that $\|p_{\mu}\|_{\infty}<\infty$, so that Step 1 applies to $\mu$ itself.

\textbf{Step 3.} The higher order formulas (\ref{representation2}) and (\ref{alternative}) are proved exactly as in \cite[Theorem 1]{BC1}, by writing $\partial_{\alpha}f=\sum_{i=1}^d(\partial_iQ_d)*(\partial_i\partial_{\alpha}f)$ and repeating the above computation.
\end{proof}

\begin{proof}[Proof of Theorem \ref{regularity}]
By (\ref{alternative}), we have
\ce
\int|\partial_{\alpha}p_{\mu_{\phi}}(x)|^pdx=\int|\partial_{\alpha}^{\mu}\phi(x)|^p|p_{\mu}(x)|^pdx\leq\|p_{\mu}\|_{\infty}^{p-1}\int|\partial_{\alpha}^{\mu}\phi(x)|^pp_{\mu}(x)dx.
\de
Then we have
\be\label{important}
\|\partial_{\alpha}p_{\mu_{\phi}}\|_{L^p}\leq\|p_{\mu}\|_{\infty}^{(p-1)/p}\|\partial_{\alpha}^{\mu}\phi\|_{L_{\mu}^p}.
\ee
For $\phi\in W_{\mu}^{m+\delta,p}$, by the discrete form of the $K$-functional (the same characterization as in Remark \ref{k}, applied to the pair $(W_{\mu}^{m,p},W_{\mu}^{m+1,p})$), for every $n$ we can find $\phi_n\in W_{\mu}^{m+1,p}$ such that
\ce
\|\phi_n-\phi\|_{W_{\mu}^{m,p}}&\leqslant& C2^{-n\delta}\|\phi\|_{W_{\mu}^{m+\delta,p}},\\
\|\phi_n\|_{W_{\mu}^{m+1,p}}&\leqslant& C2^{n(1-\delta)}\|\phi\|_{W_{\mu}^{m+\delta,p}}.
\de
Note that the alternative formula (\ref{alternative}) applies to $\phi_n$ at every order up to $m+1$: once $p_{\mu}\in L^{\infty}$ is known from Lemma \ref{boundmain}, this identity is just the integration by parts formula and only requires $\phi_n\in W_{\mu}^{m+1,p}$. Applying (\ref{important}) to $\phi_n$ and $\phi_n-\phi$ then gives
\ce
\|p_{\mu_{\phi_n}}-p_{\mu_{\phi}}\|_{W^{m,p}}
&\leqslant& C\|p_{\mu}\|_{\infty}^{(p-1)/p}2^{-n\delta}\|\phi\|_{W_{\mu}^{m+\delta,p}},\\
\|p_{\mu_{\phi_n}}\|_{W^{m+1,p}}
&\leqslant& C\|p_{\mu}\|_{\infty}^{(p-1)/p}2^{n(1-\delta)}\|\phi\|_{W_{\mu}^{m+\delta,p}}.
\de
Therefore, by the discrete $K$-functional characterization, we have
\ce
\|p_{\mu_{\phi}}\|_{W^{m+\delta',p}}\leq C\|p_{\mu}\|_{\infty}^{(p-1)/p}\|\phi\|_{W_{\mu}^{m+\delta,p}}.
\de
Combining this with (\ref{i1}), we have
\ce
\|p_{\mu_{\phi}}\|_{W^{m+\delta',p}}\leq (C(\delta,p,d,\|1\|_{W_{\mu}^{\delta,p}})\|1\|_{W_{\mu}^{\delta,p}})^{(p-1)/p}\|\phi\|_{W_{\mu}^{m+\delta,p}}.
\de
By Theorem \ref{elementary} and Lemmas \ref{pairing} and \ref{boundmain}, we have
\ce
|\partial_{\alpha}p_{\mu_{\phi}}(x)|&\leq&\sum_{i=1}^d\|\partial_iQ_d(\cdot-x)\|_{B_{\mu}^{1-\delta,p/(p-1)}}\|\partial_{(\alpha,i)}^{\mu}\phi\|_{(B_{\mu}^{1-\delta,p/(p-1)})^*}\\
&\leq&C\Gamma_p(\mu)\|\phi\|_{W_{\mu}^{m+\delta,p}}\\
&\leq&C(\delta,p,d,\|1\|_{W_{\mu}^{\delta,p}})\|1\|_{W_{\mu}^{\delta,p}}\|\phi\|_{W_{\mu}^{m+\delta,p}}.
\de
This proves the first conclusion.

By \cite[Corollary IX.13]{B} or \cite[Theorem 7]{BC1}, the second conclusion follows by a procedure similar to that in the proof of \cite[Theorem 8]{BC1}. This completes the proof.
\end{proof}

\begin{proof}[Proof of Corollary \ref{holderdensity}]
Since $\mu=\mu_{1}$, the absolute continuity of $\mu$ and the representation formula follow from Theorem \ref{elementary} applied with $\phi=1$, and the continuity of $p_{\mu}$ is also contained in that theorem. The estimate for $\|p_{\mu}\|_{\infty}$ is exactly (\ref{i1}) in Lemma \ref{boundmain} above. For the H\"{o}lder continuity, since $\delta p>d$ we can choose $\delta'\in(d/p,\delta)$ and apply Theorem \ref{regularity} (2) with $m=0$ and $\phi=1$; then $k_1=0$ and $k_2=\delta'-d/p$, and letting $\delta'\uparrow\delta$ gives the claim for every exponent $\beta<\delta-d/p$.
\end{proof}

\subsection{Sobolev spaces associated with a Wiener functional}
We denote by $\mu_F$ the law of $F$, and for a random variable $G$ we define the signed measure $\mu_{F,G}$ by $\mu_{F,G}(f):=E[f(F)G]$. For $G\in W_F^{m,p}$, we define
\ce
\partial_{\alpha}^FG=E[H_{\alpha}(F,G)|F]
\de
and
\ce
\phi_G(x)=E[G|F=x].
\de

\bd
Given a multiindex $\alpha$ and a power $p\geq 1$, we say that the integration by parts formula $IP_{\alpha,p}(F,G)$ holds if there exists a random variable $H_{\alpha}(F,G)\in\mL^p$ such that
\ce
E[\partial_{\alpha}f(F)G]=E[f(F)H_{\alpha}(F,G)],\quad\forall f\in C_c^{|\alpha|}(\mR^d).
\de
We denote by $W_F^{m,p}$ the space of random variables $G\in\mL^p$ such that $IP_{\alpha,p}(F,G)$ holds for every multiindex $\alpha$ with $|\alpha|\leq m$, equipped with the norm
\ce
\|G\|_{W_F^{m,p}}=\|G\|_{\mL^p}+\sum_{1\leq|\alpha|\leq m}E[|\partial_{\alpha}^FG|^p]^{1/p}.
\de
We define $W_F^{-k,p}=(W_F^{k,q})'$ where $1/p+1/q=1$. Let $1<p<\infty$, $\alpha\in{\mathbb R}$. We also define fractional order Sobolev spaces as follows:
\begin{eqnarray*}
{\mathcal E}_F^{\alpha,p}=(W_F^{m-1,p},W_F^{m,p})_{\theta,p}
\quad\text{where $m$ is the smallest integer larger than $\alpha$}
\end{eqnarray*}
equipped with the norm
\begin{eqnarray*}
\|u\|_{{\mathcal E}_F^{\alpha,p}}=[\int_0^{1}|\epsilon^{-\theta}K(\epsilon,u)|^p\frac{d\epsilon}{\epsilon}]^{1/p},
\end{eqnarray*}
where
\begin{eqnarray*}
K(\epsilon ,f)=\inf\{\|f_1\|_{W_F^{m-1,p}}+\epsilon \|f_2\|_{W_F^{m,p}},\; f_1+f_2=f,\; f_1\in W_F^{m-1,p},\; f_2\in W_F^{m,p}\}
\end{eqnarray*}
and $\theta:=\alpha-(m-1)$.
\ed

For $0<\delta\leqslant 1$ and $1<p<\infty$ we shall simply write
\ce
W_F^{\delta,p}:=\cE_F^{\delta,p}=(\mL^p,W_F^{1,p})_{\delta,p}.
\de

The definition of the weighted weak derivative then gives
\ce
W_F^{m,p}=\{G\in\mL^p;\phi_G\in W_{\mu_F}^{m,p}\} \quad\text{and}\quad\partial_{\alpha}^FG=(-1)^{|\alpha|}\partial_{\alpha}^{\mu_F}\phi_G(F).
\de

Taking $\phi=1$ in Theorems \ref{elementary} and \ref{regularity}, applied to the measure $\mu_F$, we obtain the following result concerning the density of $F$ itself.
\bt\label{boundness21}
Let $1<p<\infty$, $0<\delta\leqslant 1$ with $\delta p>d$ and let $1\in W_F^{\delta,p}$. Then $\mu_F(dx)=p_F(x)dx$ with
\ce
p_F(x)=-\sum_{i=1}^d\langle\partial_i^{\mu_F}1,\partial_iQ_d(\cdot-x)\rangle=\sum_{i=1}^d\langle H_i(F,1),\partial_iQ_d(F-x)\rangle,
\de
where the first pairing is the one in Lemma \ref{pairing} and the second is understood in the sense of generalized Wiener functionals (cf. Lemma \ref{pullback} below). Moreover, $p_F$ is H\"{o}lder continuous of every exponent $\beta<\delta-d/p$ and
\ce
\|p_F\|_{\infty}\leq C(\delta,p,d,\|1\|_{W_F^{\delta,p}})\|1\|_{W_F^{\delta,p}}.
\de
\et

The representation formulas and the estimates for the weighted measure $\mu_{F,G}$ now follow from Theorems \ref{elementary} and \ref{regularity}.
\bt\label{boundness22}
Let $1<p<\infty$, $0<\delta'<\delta\leqslant 1$, $m\geq0$ with $\delta p>d$ and let $1\in W_F^{\delta,p}$, $G\in W_F^{m+\delta,p}$. Then $\mu_{F,G}(dx)=p_{F,G}(x)dx$ with
\ce
p_{F,G}(x)=-\sum_{i=1}^d\langle\partial_i^{\mu_F}\phi_G,\partial_iQ_d(\cdot-x)\rangle=\sum_{i=1}^d\langle H_i(F,G),\partial_iQ_d(F-x)\rangle.
\de
Moreover, we have $p_{F,G}\in W^{m+\delta',p}$ with
\ce
\|p_{F,G}\|_{W^{m+\delta',p}}\leq (C(\delta,p,d,\|1\|_{W_F^{\delta,p}})\|1\|_{W_F^{\delta,p}})^{(p-1)/p}\|G\|_{W_F^{m+\delta,p}}
\de
and
\ce
\|p_{F,G}\|_{\infty}\leq C(\delta,p,d,\|1\|_{W_F^{\delta,p}})\|1\|_{W_F^{\delta,p}}\|G\|_{W_F^{m+\delta,p}},
\de
and for any $\alpha$ with $0\leq|\alpha|\leq m$,
\ce
\partial_{\alpha}p_{F,G}(x)=(-1)^{|\alpha|}\sum_{i=1}^d\langle\partial_{(\alpha,i)}^FG,\partial_iQ_d(\cdot-x)\rangle=(-1)^{|\alpha|}\sum_{i=1}^d\langle H_{(\alpha,i)}(F,G),\partial_iQ_d(F-x)\rangle.
\de
Furthermore, assume $m+\delta'>d/p$ and let $k_1=[m+\delta'-d/p]$ be the integer part of $m+\delta'-d/p$ and $k_2=\{m+\delta'-d/p\}$ the fractional part. When $k_2>0$, we have $p_{F,G}\in C^{k_1}(\mR^d)$; for any multiindex $\alpha$ with $0\leq|\alpha|\leq k_1-1$, $\partial_{\alpha}p_{F,G}$ is Lipschitz continuous, while for any multiindex $\alpha$ with $|\alpha|=k_1$, $\partial_{\alpha}p_{F,G}$ is H\"{o}lder continuous of exponent $k_2$. The borderline case $k_2=0$ is understood as in part (2) of Theorem \ref{regularity}. In particular, if the assumptions hold for arbitrarily large $p$, then, letting $\delta'\uparrow\delta$ and $p\uparrow\infty$ (which avoids the borderline case), we obtain that for any multiindex $\alpha$ with $0\leq|\alpha|\leq m$, $\partial_{\alpha}p_{F,G}$ is H\"{o}lder continuous of every exponent $\beta<\delta$.
\et

\section{Proof of Theorem \ref{thm:main}}
To prove Theorem \ref{thm:main}, we first establish three auxiliary lemmas. The first one is a pullback estimate, the second one is the lifting property of the scale $W_F^{\delta,p}$, and the third one bridges the Malliavin differentiability of $F$ and the hypothesis $1\in W_F^{\delta,p}$.
\bl\label{pullback}
Let $F\in\mD_1^{\infty-}(\mR^d)$, $1<q<\infty$ and $0\leqslant s\leqslant 1$. Then for every $\epsilon>0$ there exist $r>q$ and $C>0$ such that, for every $g\in W_{\mu_F}^{s,q(1+\epsilon)}$,
\ce
\|g(F)\|_{s-\epsilon,q}\leqslant C\big(1+\|\,\|DF\|_H\,\|_{r}\big)\|g\|_{W_{\mu_F}^{s,q(1+\epsilon)}}.
\de
\el
\begin{proof}
The map $g\mapsto g(F)$ is a contraction from $L_{\mu_F}^q$ into $\mL^q$. For $g\in W_{\mu_F}^{1,q}$ we have $D(g(F))=\sum_{i=1}^d\partial_ig(F)DF_i$, hence $\|D(g(F))\|_H\leqslant|\nabla g(F)|\cdot\|DF\|_H$ and, by H\"{o}lder's inequality, for every $\epsilon>0$ there is $r>q$ with
\ce
\|g(F)\|_{1,q}\leqslant C\big(1+\|\,\|DF\|_H\,\|_{r}\big)\|g\|_{W_{\mu_F}^{1,q(1+\epsilon)}}.
\de
Real interpolation between these two estimates shows that $g\mapsto g(F)$ maps $W_{\mu_F}^{s,q(1+\epsilon)}$ boundedly into $(\mL^q,\mD_1^q)_{s,q}=\cE^q_s$: the extra integrability margin $q(1+\epsilon)$ on the source side absorbs both the H\"{o}lder factor $1+\|\|DF\|_H\|_r$ and the mismatch of the Lebesgue exponents of the two endpoint spaces (see e.g. the interpolation argument in \cite[Section 1.2]{T}). Proposition \ref{two types}, with the loss of the differentiability index written as $\epsilon$, then gives $\cE^q_s\subset\mD_{s-\epsilon}^q$.
\end{proof}

\bl\label{lifting}
Let $1<p<\infty$ and $0<\delta<1$. Then there exists $C>0$ such that, for every $u\in\mL^p$ with $\partial_i^Fu\in W_F^{\delta-1,p}$, $i=1,\cdots,d$, one has $u\in W_F^{\delta,p}$ and
\ce
\|u\|_{W_F^{\delta,p}}\leqslant C\Big(\|u\|_{\mL^p}+\sum_{i=1}^d\|\partial_i^Fu\|_{W_F^{\delta-1,p}}\Big).
\de
\el
\begin{proof}
This is the lifting property of the interpolation scale $(W_F^{k,p})_{k\in\mZ}$, where $W_F^{-1,p}:=(W_F^{1,q})'$. The operators $\partial_i^F:W_F^{1,p}\to\mL^p$ and their adjoints $\partial_i^F:\mL^p\to W_F^{-1,p}$ are bounded by definition, and the scale is generated by the finite family $\{\partial_i^F\}_{i=1}^d$. Hence the retraction theorem for real interpolation (cf. \cite[Theorem 1.2.4]{T}) applies to the map $u\mapsto(u,\partial_1^Fu,\cdots,\partial_d^Fu)$ exactly as for the classical Besov scale (cf. \cite[Section 2.3.3]{T}), which gives the claim.
\end{proof}

\bl\label{bridge}
Let $0<\delta\leqslant 1$ and let $F\in\mD_{1+\delta}^{\infty-}(\mR^d)$ be nondegenerate. Then for every $p>1$ and $0<\delta'<\delta$ we have $1\in W_F^{\delta',p}$, and there exist $r>p$ and $\delta''\in(\delta',\delta)$ such that
\ce
\|1\|_{W_F^{\delta',p}}\leqslant C\Big(1+\sum_{i=1}^d\|H_i(F,1)\|_{-1+\delta'',r}\Big),
\de
where $H_i(F,1)=\delta(\gamma_F^{ij}DF_j)$ (with $j$ summed over) is the weight in the integration by parts formula.
\el
\begin{proof}
Let $\delta''\in(\delta',\delta)$ and let $q$ be the conjugate exponent of $p$. For every $g\in W_{\mu_F}^{1-\delta',q}$ with $\|g\|_{W_{\mu_F}^{1-\delta',q}}\leqslant1$, by the definition of $H_i(F,1)$ and Lemma \ref{pullback},
\ce
|E[g(F)H_i(F,1)]|\leqslant\|H_i(F,1)\|_{-1+\delta'',r}\|g(F)\|_{1-\delta'',r'}\leqslant C\|H_i(F,1)\|_{-1+\delta'',r},
\de
where $r'$ is the conjugate exponent of $r$ and $r>p$ is chosen sufficiently close to $p$ (equivalently, $r'<q$ is sufficiently close to $q$) so that Lemma \ref{pullback} applies with the integrability margin $q(1+\epsilon)=r'$. This means $\partial_i^F1\in W_F^{-1+\delta',p}$ with $\|\partial_i^F1\|_{W_F^{-1+\delta',p}}\leqslant C\|H_i(F,1)\|_{-1+\delta'',r}$ for every $i=1,\cdots,d$, and the conclusion follows from Lemma \ref{lifting} after summing over $i$.
\end{proof}

We are now in a position to complete the proof of Theorem \ref{thm:main}.
\begin{proof}[Proof of Theorem \ref{thm:main}]
(1) Under the assumptions on $F$, we have $\det(\Sigma_F)\in\mD_{\delta}^{\infty-}$. Passing from the scale $\mD_\alpha^p$ to the interpolation scale $\cE_\alpha^p$ by Proposition \ref{two types} and then applying Proposition \ref{inverse} (with the nondegeneracy assumption supplying the inverse moment estimate) gives
\ce
(\det(\Sigma_F))^{-1}\in\mD_{\delta-}^{\infty-}:=\cap_{\delta'<\delta}\mD_{\delta'}^{\infty-},
\de
which implies that $\gamma_F^{ij}\in\mD_{\delta-}^{\infty-}$. Let $f\in C_c^1(\mR^d)$. By the chain rule, we have
\ce
D(f\circ F)=\sum_{i=1}^d(\partial_if)\circ F\cdot DF_i,
\de
and then, using the summation convention,
\ce
\partial_if\circ F=\gamma_F^{ij}\cdot\<D(f\circ F),DF_j\>_H.
\de
Therefore for $G\in\mD_{\infty}^{\infty-}$, we have
\ce
E[(\partial_if)\circ F\cdot G]=E[f\circ F\cdot\delta(\gamma_F^{ij}\cdot G\cdot DF_j)]=E[f\circ F\cdot H_i(F,G)]
\de
where
\ce
H_i(F,G)=\delta(\gamma_F^{ij}\cdot G\cdot DF_j).
\de
Since $\gamma_F^{ij}\cdot DF_j\in\mD_{\delta-}^{\infty-}(H)$ and $\delta:\mD_{\alpha+1}^r(H)\rightarrow\mD_{\alpha}^r$ is continuous for every $r>1$ and $\alpha\in\mR$, we have, for every $r>1$ and $0<\delta''<\delta$,
\ce
\|H_i(F,1)\|_{-1+\delta'',r}\leq C\|\gamma_F^{ij}\cdot DF_j\|_{\delta,r}.
\de
By Lemma \ref{bridge}, for every $p>1$ and $0<\delta'<\delta$ we have $1\in W_F^{\delta',p}$ with the corresponding estimate. We now choose $p>d/\delta$ and then $\delta'\in(d/p,\delta)$, so that $\delta'p>d$; note that such a choice is possible for every dimension $d\geqslant 1$ and every $\delta>0$. Theorem \ref{boundness21} then applies and yields the absolute continuity of the law of $F$, the representation formula
\ce
p_F(x)=E[\sum_{i=1}^d\partial_iQ_d(F-x)H_i(F,1)]
\de
(understood in the sense of the pairing as in Theorem \ref{boundness21}), the H\"{o}lder continuity of $p_F$ of every exponent $\beta<\delta'-d/p$, and the estimate (\ref{w1}). Letting $p\uparrow\infty$ and then $\delta'\uparrow\delta$, we conclude that $p_F$ is H\"{o}lder continuous of every exponent $\beta<\delta$, and the first conclusion is proved.

(2) For $m\geq2$, the same argument as in part (1), using Propositions \ref{two types} and \ref{inverse}, gives $\gamma_F^{ij}\in\mD_{m-1+\delta-}^{\infty-}$. Hence $H_i(F,1)=\delta(\gamma_F^{ij}DF_j)\in\mD_{\delta-}^{\infty-}\subset\mL^p$ for every $p>1$, i.e. $1\in W_F^{1,p}$. Iterating the integration by parts as in \cite[Proposition 22]{BC1}, the higher regularity of $\gamma_F^{ij}$ just stated gives $1\in W_F^{m,p}$ for every $p>1$.

We now apply Theorem \ref{boundness22} with the integer parameter $m-1$ and the fractional parameter $\delta$. Since $\delta\leq1$, we have $1\in W_F^{m,p}\subset W_F^{(m-1)+\delta,p}$; thus, for arbitrary $p>d/\delta$ (so that $\delta p>d$) and $\delta'<\delta$, the hypotheses of Theorem \ref{boundness22} are satisfied. We obtain $p_F\in W^{m-1+\delta',p}$, the representation formulas for $\partial_{\alpha}p_F$ with $|\alpha|\leq m-1$ and the estimate (\ref{w2}). The formulas for $|\alpha|=m$ follow by one more integration by parts, understood in the generalized pairing sense as in part (1). Finally, the H\"{o}lder continuity of $\partial_{\alpha}p_F$ ($|\alpha|\leq m-1$) of every exponent $\beta<\delta$ follows from the last statement of Theorem \ref{boundness22} by taking $\delta'$ sufficiently close to $\delta$ and $p$ sufficiently large. This completes the proof.
\end{proof}

\end{document}